\documentclass[12pt]{article}
\usepackage[centertags]{amsmath}
\usepackage{amsfonts}
\usepackage{amssymb}
\usepackage{latexsym}
\usepackage{amsthm}
\usepackage{newlfont}
\usepackage{graphicx}
\usepackage{listings}
\usepackage{booktabs}
\usepackage{abstract}
\RequirePackage{srcltx}

\date{}

\def \Mm#1{\mbox{\boldmath$\scriptstyle #1$\unboldmath}}
\def \MM#1{\mbox{\boldmath$#1$\unboldmath}}

\def\baa {\begin{eqnarray*}}
\def\eaa {\end{eqnarray*}}

\def \la {\lambda}
\def \al {\alpha}

\def \ra {{\quad\Rightarrow\quad}}
\def \lr {{\quad\Leftrightarrow\quad}}

\def \k {\MM{k}}
\def \m {\MM{m}}

\def \kk {{\Mm{k}}}

\def \ub {\underbrace}

\def\R{{\mathbb R}}

\def\N{{\mathbb N}}
\def\Z{{\mathbb Z}}

\def\bk{{\mathbf k}}
\def\bx{{\mathbf x}}

\def\LL{{\cal L}}
\def\U{{\cal U}}

\def\wh{\widehat}

\def\la{\lambda}

\newtheorem{lemma}{Lemma}[section]

\newtheorem{theorem}[lemma]{Theorem}

\newtheorem{definition}[lemma]{Definition}

\theoremstyle{remark}
\newtheorem{remark}[lemma]{Remark}
\def\be  {\begin{equation}}
\def\ee  {\end{equation}}
\def\ba  {\begin{eqnarray}}
\def\ea  {\end{eqnarray}}
\def\baa {\begin{eqnarray*}}
\def\eaa {\end{eqnarray*}}
\def\bc  {}
\makeatletter
\@addtoreset{equation}{section}
\makeatother
\def\proof{\medskip\noindent{\bf Proof.} }
\def\qed{\hfill $\Box$}
\newcommand {\lb} {\label}
\newcommand {\rf[1]} {(\ref{#1})}

\begin{document}




\title{On the cardinality of lower sets and  universal discretization
	\footnote{
	\phantom{xxx} MSC 2020. {Primary: {65J05}; Secondary: 05A17, 42B05, 65D30, 41A17, 41A63 }
	\\ \phantom{xxx}
		The first named author's research was partially supported by NSERC of Canada Discovery Grant
RGPIN-2020-03909.
		The second named author's research was partially supported by NSERC of Canada Discovery Grant
RGPIN-2020-05357.
The forth
named author's research was supported by the Russian Federation Government Grant No.
14.W03.31.0031.
		The fifth  named
		author's research was partially supported by
PID2020-114948GB-I00,  2017 SGR 358,
 the CERCA Programme of the Generalitat de Catalunya,
 Severo Ochoa and Mar\'{i}a de Maeztu Program for Centers and Units of Excellence in R$\&$D (CEX2020-001084-M),
 Ministry of Education and Science of the Republic of Kazakhstan (AP09260223).
 The authors would like to thank the Isaac Newton Institute for Mathematical Sciences, Cambridge, for support and hospitality during the programme ``Approximation, sampling and compression in data science'', where work on this paper was undertaken. This program was supported by EPSRC grant no EP/R014604/1.		
 The last two authors were supported by the Simons Foundation and Pembroke  College visiting fellow grants.
}
}

\author{ F. Dai, \, A. Prymak,\, A. Shadrin, \\ V. Temlyakov, \, and  \, S. Tikhonov  }

\newcommand{\Addresses}{{
		\bigskip
		\footnotesize
		
		F.~Dai, \textsc{ Department of Mathematical and Statistical Sciences\\
			University of Alberta\\ Edmonton, Alberta T6G 2G1, Canada\\
			E-mail:} \texttt{fdai@ualberta.ca }

		\medskip
		A.~Prymak, \textsc{ Department of Mathematics\\
			University of Manitoba\\ Winnipeg, MB, R3T 2N2, Canada
			\\
			E-mail:} \texttt{Andriy.Prymak@umanitoba.ca }
		
		\medskip
		A.~Shadrin, \textsc{Department of Mathematics and Theoretical Physics\\
			University of Cambridge\\Wilberforce Road, Cambridge CB3 0WA, UK
			\\
			E-mail:} \texttt{a.shadrin@damtp.cam.ac.uk}

		\medskip
		V.N. Temlyakov, \textsc{University of South Carolina, 1523 Greene St., Columbia SC, 29208, USA; Moscow Center for Fundamental and Applied Mathematics; Steklov Institute of Mathematics; and Lomonosov Moscow State University
			\\
			E-mail:} \texttt{temlyak@math.sc.edu}
		
		\medskip
		
		S.~Tikhonov, \textsc{Centre de Recerca Matem\`{a}tica\\
			Campus de Bellaterra, Edifici C
			08193 Bellaterra (Barcelona), Spain;\\
			ICREA, Pg. Llu\'{i}s Companys 23, 08010 Barcelona, Spain,\\
			and Universitat Aut\`{o}noma de Barcelona\\
			E-mail:} \texttt{stikhonov@crm.cat}		
}}
\maketitle

\begin{abstract}
	A set $Q$ in $\mathbb{Z}_+^d$ is a lower set if $(k_1,\dots,k_d)\in Q$ implies $(l_1,\dots,l_d)\in Q$ whenever $0\le l_i\le k_i$ for all $i$. 
	We derive new and refine known results regarding the cardinality of the lower sets of size $n$ in $\mathbb{Z}_+^d$. Next we apply these results for universal discretization of the $L_2$-norm of elements from $n$-dimensional subspaces of trigonometric polynomials generated by lower sets.
\end{abstract}


\section{Introduction}


We consider here sets $Q\subset \Z^d_+$ with the following property: if $\k= (k_1,\dots,k_d)\in Q$ then the box
{$\prod_{j=1}^d [0,k_j]\cap \Z^d_+$} also belongs to $Q$.

These sets were used
in the 1980s
in the multivariate approximation of periodic functions in \cite{Tem4}, where  they  were identified as {\it $Q$ has the property $S$}.
Recently, sets with this property -- now called the {\it lower sets} or the {\it downward closed sets} --
have become important in stochastic PDEs \cite{chkifa,c-d}, in sparse multivariate
approximation  \cite{devore,cm,cmn,m0}, and in the multivariate Lagrange interpolation \cite{dyn}.

To this end, we point out that, actually, the lower sets constitute the well-known
and extensively studied part
of the number theory known as {\it integer partitions} \cite{andrews},
and the latter also got further applications (and further results)
in statistical physics \cite{g1,maslov}.




The purpose of this paper is two-fold.

Firstly, we
 derive new and refine previous results regarding the cardinality
of the lower sets in $\R^d$ of size $n$.

Secondly, we apply those results to the so-called problem of universal discretization of the $L_2$-norm
of a function from a set of $n$-dimensional trigonometric subspaces generated by lower sets.

\begin{definition}[Lower set] \rm
A {\it lower set} $Q$ in $\R^d$ is a set of {\it non-negative}
integer points $\k = (k_1,\ldots,k_d) \in\Z_+^d$
such that
$$
   \k  = (k_1,\ldots,k_d) \in Q
\ra \k' \in Q \quad\mbox{if}
\quad { \k'\in\Z_+^d \  \ \text{and}}\  \  0 \le k_i' \le k_i \quad \forall i\,.
$$
We denote by $\LL_d(n)$ the set of all lower sets in $\R^d$
of size $n$,
$$
   \LL_d(n) := \{Q \subset \Z_+^d: \mbox{$Q$ is a lower set, and $|Q| = n$} \}\,,
$$
and by $p_d(n)$ the cardinality of $\LL_d(n)$,
$$
     p_d(n) := |\LL_d(n)|\,.
$$
\end{definition}

As mentioned above, these sets are used in sparse multivariate approximation in $\R^d$ by algebraic
and trigonometric polynomials which belong to the $n$-dimensional subspaces
spanned by multinomials whose powers form a lower set (of size $n$).

Sometimes, it is more convenient to consider {\it positive} lower sets, i.e., the sets
which are obtained from the lower sets with shift by the
unit vector $\MM{1} = (1,1,\ldots,1)$.

\begin{definition}[Positive lower set] \rm
A {\it positive lower set} $P$ in $\R^d$ of size $n$ is a set of {\it positive}
integer points $\k = (k_1,\ldots,k_d) \in\N^d$
such that
$$
   \k  = (k_1,\ldots,k_d) \in P
\ra \k' \in P \quad\mbox{if}
\quad  { \k'\in \N^d \  \ \text{and}}\  \  1 \le k_i' \le k_i \quad \forall i
$$
and $|P| = n$.

In other words, if an integer point $\k \in \N^d$ belongs to a positive lower set $P$,
then the hypercube $H_{\kk} = \prod_{i=1}^d [0,k_i] \cap \N^d$
also belongs to $P$.
Thus, any positive lower set $P$ is a union of some
finite number of hypercubes $H_\kk$ in $\N^d$.
\end{definition}

With a slightly different definition, the positive lower sets in $\N^d$ are well-known
in number theory as integer partitions, with a vast literature
on this subject (though restricted mainly to the cases $d=2$ and $d=3$), see \cite{andrews}.

\begin{definition}[Integer partition in $\R^d$] \rm
An {\it integer partition} $P$ of $n \in \N$ in $\R^d$ is any representation of $n$ in the form
\be \lb{np}
     n = \sum_{\kk \in \N^{d-1}} n_{\kk}, \quad \mbox{where} \;\; n_\kk \in \N
     \; \mbox{and}\; n_{\kk'} \ge n_{\kk} \;\;\mbox{if}\;\; 1 \le k_i' \le k_i \quad \forall i, \quad 1\le i \le d-1\,.
\ee
\end{definition}

We may visualize an integer partition $P$ of $n \in \N$ in $\R^d$ as a union of stacks of $d$-dimensional unit cubes,
where each stack consists of $n_\kk$ cubes one over the other, with the lowest one having the $(d-1)$-dimensional
base centered at $\k - \MM{\frac{1}{2}}$.

Now, with any $d$-dimensional partition of $n$ (that satisfies \rf[np]),
$$
     n = \sum_{\kk \in \N^{d-1}} n_\kk\,,
$$
we can identify the set of points $\wh \k \in \N^d$ by the rule
$$
      \Big(n_\kk \in \N,\; \k = (k_1,\ldots,k_{d-1}) \in \N^{d-1} \Big)
\quad\to\quad \wh\k = (k_1,\ldots,k_{d-1}, 1 \le k_d \le n_\kk) \in \N^d\,,
$$
and clearly this gives a one-to-one correspondence between integer partitions of $n$ in $\R^d$
and positive lower sets of size $n$ in $\R^d$.

Besides number theory, integer partitions also prominently appear in statistical mechanics (see \cite{g3,maslov}),
and in studying the so-called lattice animals (also known as polyominoes) \cite{bpa}. In fact, all known
results about the cardinality $p_d(n)$ of integer partitions for high dimensions $d > 3$
were obtained by those working in these areas.

We also mention that, from the definition, it follows that an integer partition
is a monotonely decreasing integer-valued function
of $d$ integer arguments, and this connection was studied in \cite{bolour}.

\medskip
{\bf 1.1. Integer partitions for  $d \le 4$.}

1) A {\it linear or two-dimensional partition} of a positive integer $n$ is given by a formula
$$
   n = n_1 + n_2 + \cdots + n_k\,, \qquad \mbox{where $n_i \ge n_{i+1}$}\,.
$$
The asymptotic number of all such partitions $p_2(n)$ as $n \to \infty$ is given by the celebrated Hardy-Ramanujan
result \cite{hr} (1918),
\be \lb{d=2}
   p_2(n) \sim \frac{1}{4n\sqrt{3}} \exp (\alpha_2 \sqrt{n}), \qquad
   \alpha_2 = 2 \zeta(2)^{1/2} = \pi \sqrt{\frac{2}{3}}\,,
\ee
so that
\be \lb{a_2}
      \lim_{n\to \infty} \frac{\ln p_2(n)}{n^{1/2} } = \alpha_2 \approx 2.565099.
\ee
Here, $\zeta(s) := \sum_{n=1}^\infty \frac{1}{n^s}$ is the Riemann zeta-function.
Actually, the asymptotic bound provides an upper bound for all $n$ \cite{p}, namely
\be \lb{upper2}
    p_2(n) \le e^{\alpha_2 \sqrt{n}}\,, \quad n \in \N\,.
\ee

\medskip
2) Similarly, a {\it plane or three-dimensional partition} of $n \in \N$ is given by a formula
$$
   n = \sum_{i,j} n_{i,j}, \quad
  \mbox{where}  \quad n_{i,j} \ge n_{i+1,j}, \; n_{i,j} \ge n_{i,j+1}\,.
$$
In this case the asymptotics (as $n \to \infty$) was found by Wright \cite{w} (1931)
\be \lb{d=3}
        p_3(n) \sim c_3 {n}^{-25/36} \exp \big(\alpha_3 n^{2/3} \big), \qquad
        \alpha_3 = \frac{3}{2} \Big(2\zeta(3)\Big)^{1/3},
\ee
with explicitly given constant $c_3$ and the following numerical value for $\alpha_3$
\be \lb{a_3}
      \lim_{n\to \infty} \frac{ \ln p_3(n) }{ n^{2/3} } = \alpha_3 \approx 2.00945\,.
\ee

3) Finally, a {\it solid or four-dimensional partition} of $n$ is given by a formula
$$
      n = \sum_{i,j,k} n_{i,j,k},
\quad \mbox{where} \quad n_{i,j,k} \ge \max \Big( n_{i+1,j,k}, n_{i,j+1,k}, n_{i,j,k+1} \Big).
$$
In that case, the exact asymptotics as $n \to \infty$ is unknown. There was a conjecture \cite[Conjecture 3.1]{g3} that,
similarly to the cases $d=2,3$ in \rf[d=2] and \rf[d=3], it should be given for all $d \ge 2$ by
the following expression
\be \lb{d-lim}
      \alpha_d := \lim_{n\to \infty} \frac{\ln p_d(n)} { n^{1-1/d} }
\;\stackrel{?}{=}\; \frac{d}{d-1} \Big[(d-1)\zeta(d)\Big]^{1/d} =: \rho_d\,.
\ee
However, that cannot be true for large $d$, since $\rho_d \to 1$ as $d \to \infty$,
whereas we show in this paper (Theorem \ref{thm2}) that
$\alpha_d > \frac{d}{(d!)^{1/d}}\ln 2 \to e \ln 2 = 1.8841$.
In fact, from our lower estimate, it follows that
$$
    \alpha_d > \rho_d \quad \mbox{for} \; d \ge 8.
$$
Also, for $d=4$,
the Monte-Carlo simulations \cite{g2} produced for $\alpha_4$
a different 
value,
\be \lb{a_4}
      \alpha_4 = \lim_{n\to \infty} \frac{ \ln p_4(n) }{ n^{3/4} }  \approx 1.822,
      \quad \rho_4 \approx 1.78982.
\ee

\medskip
{\bf 1.2. Integer partitions for $d > 3$.}
For $d > 3$, Bhatia, Prasad and Arora \cite{bpa} (1997) proved that,
for a fixed $d$ and for sufficiently large $n$, we have
\be \lb{as}
   C_1(d) \le \frac {\ln p_d(n)}{n^{1-1/d}} \le C_2(d)\,, \qquad n > n_d\,,
\ee
with some constants $C_1(d),C_2(d)$ that depend on $d$. However,
they did not write down explicitly the nature of this dependence, and we will comment on the bounds that can be obtained from their work in Remarks~\ref{rem:bhatiaproof} and~\ref{rem:lower}.

Note also that some papers (e.g. \cite{g3}) claim that in \cite{bpa} existence of the limit
$\lim\limits_{n\to \infty} \frac {\ln p_d(n)}{n^{1-1/d}}$ was proved,
but that was not the case. In this paper, we use the same limit notation but we mean by that ``$\limsup$'' for the upper bounds, and ``$\liminf$'' for the lower bounds.

In \cite{cmn}, Cohen, Migliorati and Nobile, while studying lower sets, proved that, uniformly in $d$ and $n$,
we have
\be \lb{cohen}
       p_d(n) \le 2^{dn}, \qquad p_d(n) \le d^{n-1} (n-1)!
\ee
whichever bound is better.

We mention that, with exception of the upper bound \rf[upper2] for $d=2$, estimates \rf[cohen] were the only
explicit upper bounds for $p_d(n)$ obtained so far -- all previous estimates for the cardinality of integer partitions
for $d \ge 3$ were asymptotic, with dimension $d$ being fixed and the value of $n$  (which is being partitioned) going to infinity.

Also, there were several results concerning calculation of the exact values of $p_d(n)$ for relatively small $n$
(those values grow very fast with $n$).
The basis for these calculations is the Knuth algorithm \cite{k}, with several important contributions from
Govindarajan et al. \cite{g3, g2, g1}. In the arXiv version of \cite{g1}, a table with the values of $p_d(n)$ is given
for all $d \le 11$ and $n \le 24$.

\medskip
{\bf 1.3. New results.}
In this paper, we improve the upper estimates \rf[cohen],
and we also find explicit expressions for constants $C_1(d), C_2(d)$ in \rf[as] such that \rf[as]
is valid not just for large $n > n_d$ with unspecified $n_d$, but for all $n \in \N$.

Our first result provides lower and upper bounds which are uniform in $d$ and $n$,
and which improve the bounds \rf[cohen] by Cohen-Migliorati-Nobile.

\begin{theorem} \lb{thm1}
For any $d \ge 2$  and $n \in \N$, we have
$$
     \frac{1}{(n-1)!}\, d^{n-1} < p_d(n) \le d^{n-1}\,.
%
$$
\end{theorem}

Thus, for a fixed $n \in \N$, we have
$$
p_d(n) \asymp d^{n-1}, \quad \text{as} \quad d \to \infty.
$$
Actually, it follows from results of Atkin et al. \cite{At} that 
$$
p_d(n) =\frac{d^{n-1}}{(n-1)!}+O(d^{n-2}), \quad \text{as} \quad d \to \infty.
$$

Our second result is a refinement of \rf[as]. We follow the same scheme as suggested by Bhatia, Prasad and Arora \cite{bpa}
but we make substantial efforts to find explicit expressions for the constants involved, as well as to get the result
valid for all $n\in\N$.

\begin{theorem}  \lb{thm2}
For any $d \ge 2$  and $n \ge 2$, we have
$$
     C'_{d,n} \le \frac{\ln p_d(n)}{n^{1-1/d}} \le C_d,
$$
where
\ba
      && C'_{d,n} = \Big(1 - \min \{{\textstyle \frac{d}{d+1}}, e n^{-1/d}\} \Big)^2 \la_d \ln 2\,,
      \quad \la_d := \frac{d}{(d!)^{1/d}}\,, \nonumber \\
      && C_d =  \alpha_2 d^{\ln d}, \quad \alpha_2 = \pi \sqrt{\frac{2}{3}}\,. \lb{C}
\ea
In particular, for $\alpha_d :=\displaystyle \lim_{n\to \infty} \tfrac{p_d(n)}{n^{1-1/d}}$, 
we have the following lower and upper bounds
\be
\label{eqn:alpha_growth}
\frac{d}{(d!)^{1/d}} \ln 2 < \alpha_d < \alpha_2 d^{\ln d}.
\ee
\end{theorem}


We note that, from our lower bound, it follows that 
$$
\lim_{d\to\infty}  \alpha_d >  e \ln 2 \approx 1.8841 > 1.822 \approx \alpha_4\,,
$$
and that means that either the value of $\alpha_4$ in \rf[a_4] is wrong, or if it is correct, then the numbers
$$
\alpha_d = \lim\limits_{n\to \infty} \frac{ \ln p_d(n) }{ n^{1-1/d} }
$$
do not decrease monotonically in $d$
as the values $\alpha_2 > \alpha_3 > \alpha_4$ in \rf[a_2], \rf[a_3], \rf[a_4] may suggest. 

We do not know if $\alpha_d$'s actually increase with large $d$ as the upper estimate in \rf[eqn:alpha_growth] implies,
but even if they do, we believe that it happens at a much slower rate than given in \rf[eqn:alpha_growth].

Theorem \ref{thm1} is proved in Section 2, and Theorem \ref{thm2} in Sections 3-4.
In Section 5, we apply these results to derive a bound
for the cardinality of the so-called universal discretization point set $(\xi_i)_{i=1}^m$
for collection of trigonometric
polynomials in $d$ variables whose harmonics form a lower set of size $n$.


\section{Proof of Theorem \ref{thm1}}


We split Theorem \ref{thm1} in two parts -- the upper  and the lower estimates.

\begin{theorem}
We have
\be \lb{thm3}
   p_d(n) \le d^{n-1}\,.
\ee
\end{theorem}

\proof
The proof is by induction on $d$.

1) For $d=2$, we may partition $n$ as follows. Firstly we write
$$
   n = \underbrace{1 + 1 + 1 \cdots + 1 + 1}_{\mbox{$n$ times}}
$$
and then we form a partition of $n$ in $k$ summands as
$n = n_1 + n_2 + \cdots + n_k$ by putting $k-1$ splits in the sum above
$$
   n = \ub{1 + \cdots + 1}_{n_1} \stackrel{1}{\oplus}
       \ub{1 + \cdots + 1}_{n_2} \stackrel{2}{\oplus} \cdots
       \ub{1 + \cdots + 1}_{n_{k-1}} \stackrel{k-1}{\oplus}
       \ub{1 + \cdots + 1}_{n_k}
$$
For each $k$ we have ${n-1 \choose k-1}$ places where to put the split in,
so the total number of such (non-ordered) partitions is clearly
$$
           \sum_{k=1}^{n} {n-1 \choose k-1} = 2^{n-1}.
$$
Therefore, since integer partitions are the ordered ones, we have
$$
    p_2(n) \le 2^{n-1}\,.
$$

2) Assuming induction hypothesis to be true for $d$, we make the cuts
of each $(d+1)$-dimensional partition along the hyperplanes
perpendicular to the $x_1$-axis, say.
Let there be $k-1$ such cuts, with each $d$-dimensional part
contaning $n_i$ points, with $\sum_{i=1}^k n_i = n$.
Then the total number of $(d+1)$-dimensional partitions can be estimated as follows,
\baa
      p_{d+1}(n)
&\le& \sum_{k=1}^{n} {n-1 \choose k-1} \prod_{i=1}^k p_d(n_i) \\
&\stackrel{\rf[thm3]}{\le}& \sum_{k=1}^{n} {n-1 \choose k-1} \prod_{i=1}^k d^{n_i-1} \\
& = & \sum_{k=1}^{n} {n-1 \choose k-1} d^{n-k} \\
& = & \sum_{k=1}^{n} {n-1 \choose k-1} d^{(n-1)-(k-1)} \cdot 1^{k-1} \\
& = & (d+1)^{n-1}\,.
\eaa
This proves the upper bound for $d+1$, hence the theorem.
\qed

\medskip
For the lower bound, we will use the following well-known lemma.

\begin{lemma} \lb{le2}
Let $w_d(m)$ be the number of ways of distributing $m$ balls between $d$ ordered boxes. Then
\be \lb{w}
    w_d(m) = {m + (d-1) \choose d-1}.
\ee
\end{lemma}

\begin{theorem}
We have
$$
    p_d(n)
\ge { (d-1)+ (n-1) \choose n-1}
  = \frac{1}{(n-1)!}\, d (d+1) \cdots (d+n-2) \ge \frac{1}{(n-1)!}\, d^{n-1}\,.
$$
\end{theorem}

\proof
We obtain the lower bound as the cardinality of the subset of lower sets where all the points lie on coordinate axes.
This means that we put the first point at the origin, and then distribute the remaining $n-1$ points along $d$ axes.
This gives the required bound using the lemma above.
\qed


\section{Proof of Theorem \ref{thm2}: the lower bound}


\begin{theorem}
For any $d \ge 2$  and $n \ge 2$, we have
$$
      \frac{ \ln p_d(n) }{ n^{1-1/d}}
\ge C_d' = \Big(1- \min \{ {\textstyle\frac{d}{d+1}} , en^{-1/d}\} \Big)^2 \la_d \ln 2,
$$
where
$$
      \la_d := \frac{d}{(d!)^{1/d}} \to e\;\; (d \to \infty)\,.
$$
\end{theorem}

\proof
1) Given $m\in \Z_+$, consider the sets
\be \lb{A_m}
    A_m := \{\k \in \Z_+^d: k_1 + k_2 + \cdots + k_d = m\}\,,
\ee
and
$$
    B_m := \{\k \in \Z_+^d: k_1 + k_2 + \cdots + k_d \le m\}\,.
$$
So, for $d=2$, $B_m$ is a staircase, for $d=3$ it is  a pyramid, etc.
Clearly, $B_m$ is a lower set, and $A_m$ may be viewed as a set of ``corners" (or vertices)
of $B_m$.

\smallskip
2) Let $a_m$ and $b_m$ be the cardinality of $A_m$ and $B_m$, respectively,
$$
     a_m := |A_m|, \qquad b_m := |B_m|\,.
$$
If we remove any number of corner points (vertices)
from $B_m$, 
then the remaining set of points will form a lower set as well,
so collection of all such sets (where any corner is either ``in" or ``out") has the cardinality $2^{a_m}$.

\smallskip
3) Next, given $n \ge 2$, we choose $m\in \Z_+$ (note that $b_0=1$) such that
$$
   b_m < n \le b_{m+1}\,,
$$
and with such an $n$ we take all the lower sets constructed in step 2 above
which contain the point (corner) $(m,0,\ldots,0)$ (i.e. this corner is always ``in").
All these lower sets are subsets of $B_m$, with $< n$ points in any of them,
and we put all the remaining points  along $x_1$-axis above the vertex $(m,0,\ldots,0)$,
so that the total cardinality of any such set is exactly $n$.
This gives
$$
       p_d(n) \ge 2^{a_m - 1} = \frac{1}{2}\, 2^{a_m}\,.
$$
We can choose any of the $d$ axes where to keep the corner and
put remaining points on, hence
\be \lb{a}
       p_d(n) \ge \frac{d}{2}\, 2^{a_m} \ge 2^{a_m}\,.
\ee

It remains to find relation between $a_m$ and $n$.

\smallskip
4) From the definition \rf[A_m] of $A_m$, we see that $a_m := |A_m|$ is the number of ways of
putting $m$ balls into $d$ boxes, hence by Lemma \ref{le2}
$$
   a_m = {m+(d-1) \choose (d-1)}\,,
$$
and respectively,
$$
  b_m = \sum_{s=0}^m a_s
= \sum_{s=0}^m {s+(d-1) \choose (d-1)}
= {m+d \choose d}
= \frac{(m+1)\cdots(m+d)}{d!}\,.
$$
We also have
\be \lb{b}
   b_m = a_m \frac{m+d}{d}\,, \qquad b_m d! \le (m+d)^d.
\ee

\smallskip
5) Further,
$$
  b_m^{1-1/d}
= \frac{b_m}{b_m^{1/d}}
= a_m \frac{m+d}{d} \frac{(d!)^{1/d}}{[(m+1)\cdots(m+d)]^{1/d}}
\le a_m \frac{(d!)^{1/d}}{d} \frac{m+d}{m+1}\,,
$$
hence, with $\la_d= \frac{d}{(d!)^{1/d}}$,
$$
   a_m \ge \la_d\, \frac{m+1}{m+d}\, b_m^{1-1/d}\,. 
$$
Since $b_m = b_{m+1}\frac{m+1}{m+d+1}$, we get
$$
   a_m \ge \la_d b_{m+1}^{1-1/d}\,
   \frac{m+1}{m+d}\Big(\frac{m+1}{m+d+1}\Big)^{1-1/d}
   > \la_d b_{m+1}^{1-1/d} \Big(1 - \frac{d}{m+d+1} \Big)^2\,.
$$

6) To estimate the last factor, we use
$$
     b_{m+1} \ge n \ra m+d+1 \stackrel{\rf[b]}{\ge} [b_{m+1} d!]^{1/d} \ge (nd!)^{1/d}\,,
$$
hence
$$
    \frac{d}{m+d+1}
    \le \min \left\{ \frac{d}{d+1}, \frac{d}{(d!)^{1/d} n^{1/d}}  \right\}
    \le \min \left\{ \frac{d}{d+1}, \frac{e}{n^{1/d}}  \right\},
$$
and we obtain
$$
     a_m > \la_d n^{1-1/d} \Big(1 - \min\{ {\textstyle \frac{d}{d+1}}, e n^{-1/d} \}\Big)^2.
$$
It follows now from \rf[a] that
$$
     \frac{\ln p_d(n)}{n^{1-1/d}}
     >  \Big(1 - \min \{{\textstyle \frac{d}{d+1}}, en^{-1/d} \}\Big)^2 \la_d \ln 2\,,
$$
and that proves the theorem.
\qed

\begin{remark}\label{rem:lower}
Let us compare our lower bound for $p_d(n)$  with the one that can be derived 
from the proof of Bhatia et al. \cite{bpa}. With $\alpha_d=\displaystyle\lim_{n\to\infty}\tfrac{\ln p_d(n)}{n^{1-1/d}}$, the proof in \cite{bpa} implies that for $d\to\infty$ we have
	$
	\alpha_d>\left(\frac2{d-1}\right)^{1-1/d}\ln 2\to0,
	$
	whereas for our bound $\alpha_d>\frac{d}{(d!)^{1/d}}\ln 2\to e\ln 2\approx 1.8841$. Not only ours is a better bound, it also shows that conjecture \rf[d-lim] is wrong, and that monotonic decrease of $\alpha_d$ is questionable (as pointed out in Introduction).
\end{remark}


\section{Proof of Theorem \ref{thm2}: upper estimate}


\begin{theorem} \lb{thm4}
For all $d \ge 2$ and all $n \in \N$, we have
\be \lb{up1}
   p_d(n) \le (\gamma_d)^{n^{1-1/d}}\,, \qquad \gamma_d = (\beta_2)^{ d^{\ln d}},
\ee
where
$$
       \beta_2 = e^{\al_2} = e^{\pi \sqrt{2/3}} \approx 13.0019
$$
is the Hardy-Ramanujan constant in \rf[d=2] for $d=2$.
\end{theorem}

\proof
The proof is by induction on $d$, with the statement being true for $d=2$, since by \rf[upper2]
$$
       p_2(n) \le (\beta_2)^{\sqrt{n}}, \quad \mbox{and} \quad \beta_2 < \gamma_2 = (\beta_2)^{2^{\ln 2}}\,.
$$
So, given $d \ge 3$, assume that \rf[up1] holds for  the value $d-1$, namely for all $n\in \N$ we have
\be \lb{up2}
   p_{d-1}(n) \le (\gamma_{d-1})^{n^{1-1/(d-1)}}\,, \qquad \gamma_{d-1} = (\beta_2)^{ (d-1)^{\ln (d-1)}}.
\ee

1) Firstly, we claim that
\be \lb{sigma}
      \mbox{(i)} \quad d^{\sigma_1 d} \le \gamma_{d-1}\,, \quad \sigma_1 = 1.25, \qquad\qquad
      \mbox{(ii)}\quad d^{\sigma_2 d} \le \gamma_d\,, \quad \sigma_2 = 2.6\,.
\ee
Indeed, let us rewrite the first inequality (i)  in the equivalent form
$$
      d^{\sigma_1 d} \le \gamma_{d-1}
\lr   d^{\sigma_1 d} \le e^{\alpha_2 (d-1)^{\ln (d-1)}}
\lr  \sigma_1 \le \frac{\alpha_2 (d-1)^{\ln (d-1)}} {d \ln d} =: f_1(d)\,.
$$
It is easy to verify that $f_1' > 0$ for $d \ge 3$, hence with $\sigma_1 < f_1(3) \approx 1.2583$,
we have $\sigma_1 < f_1(d)$, and part (i) of \rf[sigma] follows.

Similarly, for the second inequality (ii) in \rf[sigma],
we need $\sigma_2 \le \frac{\alpha_2 d^{\ln d}} {d \ln d} =: f_2(d)$, where also $f_2' > 0$ for $d \ge 3$,
so we can take $\sigma_2 < f_2(3) \approx 2.6020$.

\medskip
2) From estimate \rf[thm3], and from \rf[sigma], it follows that if $n < (\sigma_2 d)^d$, then
$$
      p_d(n) < d^n = d^{ n^{1/d} n^{1-1/d} } <  d^{\sigma_2 d \cdot n^{1-1/d} } < \gamma_d^{n^{1-1/d}},
$$
hence \rf[up1] is true.
Therefore, from now on we assume that
\be \lb{n}
      n \ge (\sigma_2 d)^d\,.
\ee

\smallskip
3) For a positive lower set $A$ of size $n$ let $m_0$ be the largest integer
such that $A$ contains the point $\m_0 = (m_0,m_0,\ldots,m_0) \in \N^d$.
Since $|A| = n$, we clearly have
$$
       m_0 \le m_1 := \lfloor n^{1/d} \rfloor.
$$
Then $A$ is contained in the union of $m_2 := d m_1$ hyperplanes
$$
   \{x_1 = 1\}, \ldots, \{x_d = 1\}, \quad
   \{x_1 = 2\}, \ldots, \{x_d = 2\}, \quad\cdots\quad
   \{x_1 = m_1\}, \ldots, \{x_d = m_1\}\,.
$$
Next, we start to cut $A$ into $(d-1)$-dimensional slices (orthogonal to one of the coordinate axes)
choosing each time the hyperplane at the direction where it contains the largest
number of remaining points from $A$.
Let this number be $n_i$ at the $i$-th step, so that
$\sum_{i=1}^{m_2} n_i = n$ and $n_i \ge n_{i+1}$.

The number $w$ of ways (directions) in which such $m_2$ hyperplanes can appear in our algorithm is trivially bounded by
$$
      w \le d^{m_2} = d^{dm_1}\,,
$$
as at each of possible $m_2$ steps we can make a cut in no more than $d$ directions.

At the $i$-th step, the $(d-1)$-dimensional slice contains $n_i$ points which
form a positive lower set in $\R^{d-1}$,
so we have $p_{d-1}(n_i)$ possibilities in this slice, hence the total number of possible
lower sets in $d$ dimensions for each particular arrangement of slices
will be bounded by $\prod_{i=1}^{m_2} p_{d-1}(n_i)$. (Many of the resulting sets in $\R^d$
will not be lower sets, so this is a rather rough upper bound.)


Thus, we obtain
\be \lb{upper}
      p_d(n) \le d^{m_2} \sum\nolimits' \prod_{i=1}^{m_2} p_{d-1}(n_i)\,.
\ee
Now, we estimate (and explain) all the factors involved.

\smallskip
4) For the first factor, we have
\be \lb{d}
     d^{m_2} = d^{dm_1} \le d^{dn^{1/d}}\,.
\ee

5) For the sum, by construction, we have $n_i \ge n_{i+1}$ and
$\sum_{i=1}^{m_2} n_i = n$, i.e., the sum in \rf[upper] extends
over all values of $n_i$ which form an integer (linear) partition of $n$ in no more than $m_2$ summands.
By the upper bound \rf[upper2] for linear partitions (with any number of summands), we have
\be \lb{sum}
   \sum\nolimits' 1 \le (\beta_2)^{n^{1/2}}\,.
\ee

\smallskip
6) For the product in \rf[upper], using the induction hypothesis \rf[up2], we obtain
$$
    \prod_{i=1}^{m_2} p_{d-1}(n_i)
\le \prod_{i=1}^{m_2} (\gamma_{d-1})^{n_i^{1-1/(d-1)}}
 =  \gamma_{d-1}^s, \qquad s = \sum_{i=1}^{m_2} n_i^{1-1/(d-1)}\,.
$$
To estimate $s$, we apply the H\"older inequality
$$
      \sum_{i=1}^N |a_i| \le N^{1/p} \Big\{\sum_{i=1}^N |a_i|^q\Big\}^{1/q}\,, \qquad 1/p + 1/q = 1\,,
$$
with $\frac{1}{p} = \frac{1}{d-1}$, $\frac{1}{q} = 1 - \frac{1}{d-1}$, $a_i = n_i^{1/q}$ and $N = m_2 \le dn^{1/d}$.
This gives
\baa
       s
= \sum_{i=1}^{m_2} n_i^{1- 1/(d-1)}
&\le& (dn^{1/d})^{1/(d-1)} \Big(\sum_{i=1}^{m_2} n_i\Big)^{1-1/(d-1)} \\
&=& (dn^{1/d})^{1/(d-1)} \cdot n^{1-1/(d-1)} \\
&=& d^{1/(d-1)} n^{1-1/d} \,,
\eaa
i.e.,
\be \lb{prod}
      \prod_{i=1}^{m_2} p_{d-1}(n_i) \le (\gamma_{d-1})^{ d^{1/(d-1)} n^{1-1/d} }.
\ee

\smallskip
7) So, from \rf[upper], using \rf[d]-\rf[prod], we obtain
\be\lb{eqn:ihformula}
     p_d(n)
\le (\beta_2)^{n^{1/2}} d^{dn^{1/d}} (\gamma_{d-1})^{ d^{1/(d-1)} n^{1-1/d}}
 =: (\la_d)^{n^{1-1/d}}\,,
\ee
where
\be \lb{la}
     \la_d = (\beta_2)^{1/n^{1/2-1/d}} \cdot d^{d/n^{1-2/d}} \cdot (\gamma_{d-1})^{d^{1/(d-1)}} \,.
\ee

8) Let us bring all the factors in \rf[la] to the powers of $\gamma_{d-1}$. By \rf[up2] and \rf[sigma], we have
$$
     \beta_2 = (\gamma_{d-1})^{1/ (d-1)^{\ln (d-1)}}\,, \qquad
     d^d < (\gamma_{d-1})^{1/\sigma_1},
$$
and from \rf[n]
$$
      n > (\sigma_2 d)^d
\quad\ra\quad \frac{1}{n^{1/2-1/d}} < \frac{1}{(\sigma_2 d)^{d/2-1}}, \qquad
     \frac{1}{n^{1-2/d}} < \frac{1}{(\sigma_2 d)^{d-2}}\,.
$$
It follows from \rf[la] that
$$
     p_d(n)  < \Big(\gamma_{d-1}^{r(d)}\Big)^{n^{1-1/d}}
$$
where
\be \lb{r0}
     r(d) = \frac{1}{ (\sigma_2 d)^{d/2-1} (d-1)^{\ln (d-1)} } + \frac{1}{\sigma_1 (\sigma_2 d)^{d-2}} + d^{1/(d-1)}.
\ee
Thus, to prove that $p_d(n) < \gamma_d^{n^{1-1/d}}$, we need to check the inequality
\be \lb{final}
      (\gamma_{d-1})^{r(d)} < \gamma_d
 \lr  (d-1)^{\ln (d-1)} r(d) < d^{\ln d}
\lr    r(d) <  \frac{d^{\ln d}}  {(d-1)^{\ln (d-1)} }\,.
\ee

9) We remind that in \rf[r0], we have the values  $\sigma_1 = 1.25$, $\sigma_2 = 2.6$.

If $d=3$, then
$$
     r(3) = \frac{1}{(2.6\cdot 3)^{1/2} 2^{\ln 2}} + \frac{1}{1.25 \cdot 2.6 \cdot 3} + 3^{1/2} < 2.057,
$$
and
$$
     \frac{3^{\ln 3}}{2^{\ln 2}} > 2.06,
$$
so \rf[final] is true.

If $d \ge 4$, then $\ln (d-1) > 1$, hence $(d-1)^{\ln (d-1) } > d-1 \ge \frac{3}{4} d$, and for the first term in \rf[r0] we have
$$
      \frac{1}{(\sigma_2d)^{d/2-1} (d-1)^{\ln (d-1)} }
 <  \frac{1}{ 2.6 d \cdot 0.75 d} < \frac{0.6}{d^2}\,,
$$
whereas for the second,
$$
         \frac{1} {\sigma_1 (\sigma_2 d)^{d-2}} < \frac{1}{1.25 (2.6 d)^2} < \frac{0.2}{d^2}\,,
$$
so from \rf[r0] we obtain
\be \lb{r}
     r(d) < \frac{1}{d^2} + d^{1/(d-1)}\,.
\ee
%
%

From the inequality
$$
    (x+1)^{1/m} - x^{1/m} > \big(t^{1/m}\big)'\big|_{t=x+1} = \frac{1}{m(x+1)^{1-1/m}} > \frac{1}{m(x+1)}\,,
$$
it follows that
$$
    (d+1)^{1/(d-1)} - d^{1/(d-1)}  >  \frac{1}{(d-1)(d+1)} > \frac{1}{d^2}\,,
$$
therefore, from \rf[r],
\be\lb{eqn:newbd}
      r(d) < \frac{1}{d^2} + d^{1/(d-1)} < (d+1)^{1/(d-1)}\,.
\ee

10) Thus, to prove \rf[final] for $d \ge 4$, we need to prove that
$$
         (d+1)^{1/(d-1)} <  \frac{d^{\ln d}} {(d-1)^{\ln (d-1)}}\,.
$$
Taking natural logarithm of both sides, we obtain
$$
      \frac{\ln(d+1)}{d-1} < \ln^2 d - \ln^2 (d-1)\,,
$$
and this is true since
$$
     \ln^2 d - \ln^2 (d-1)  >  (\ln^2 x)'\big|_{x=d} = \frac{2 \ln d}{d} \,,
$$
and $\frac{2 \ln d}{d} > \frac{\ln(d+1)}{d-1}$ for $d \ge 4$.
\qed

\begin{remark}\label{rem:bhatiaproof}
Now we compare our upper bound for $p_d(n)$  with the bound that can be obtained following the proof of Bhatia et al. \cite{bpa}. Our proof is following the arguments in \cite{bpa} up to 
Eq \rf[eqn:ihformula], and they obtained an asymptotic bound for $p_d(n)$, the constant $C_2(d)$ in \rf[as], by iterating \rf[eqn:ihformula]. 

Such an iteration of \rf[eqn:ihformula] (with respect to $d$), besides the highest exponent of order
$n^{1-1/d}$, produces also a growing number of other exponents of lower order which 
look as follows
\baa
&& d = 2, \qquad c_{2,1}^{n^{1/2}};  \\
&& d = 3, \qquad c_{3,1}^{n^{1/3}} c_{3,2}^{n^{1/2}} c_{3,3}^{n^{2/3}}; \\
&& d = 4, \qquad c_{4,1}^{n^{1/4}} c_{4,2}^{n^{1/2}} c_{4,3}^{n^{5/8}}  c_{4,4}^{n^{3/4}}; \\
&& d = 5, \qquad c_{5,1}^{n^{1/5}} c_{5,2}^{n^{2/5}} c_{5,3}^{n^{1/2}}  c_{5,4}^{n^{3/5}} 
c_{5,5}^{n^{7/10}}  c_{5,6}^{n^{4/5}};\\
&& d = 6, \qquad c_{6,1}^{n^{1/6}} c_{6,2}^{n^{1/3}} c_{6,3}^{n^{1/2}}  c_{6,4}^{n^{7/12}} 
c_{6,5}^{n^{2/3}}  c_{6,6}^{n^{3/4}} c_{6,7}^{n^{5/6}}\,.
\eaa                           
Here all the constants $c_{d,k}$ will grow with $d$. For the constant with the 
highest exponent $n^{1-1/d}$,  
one can see from \rf[eqn:ihformula] that $\alpha_d=\displaystyle\lim_{n\to\infty}\tfrac{\ln p_d(n)}{n^{1-1/d}}$ will grow as follows
$$
\alpha_d \le d^{1/(d-1)} \alpha_{d-1},
$$
and this results in the estimate 
$$
\alpha_d \le \alpha_2 \prod_{k=3}^d k^{1/(k-1)}, \qquad d \ge 3.
$$
It follows then that
\be\lb{eqn:1}
\ln \alpha_d   \le  \ln \alpha_2 + \sum_{k=3}^d \frac{\ln k}{k-1},
\ee
and for the sum one can get the following bound
$$
\sum_{k=3}^d \frac{\ln k}{k-1}  \le  \int_2^d \frac{\ln x}{x}\,dx + c_1
\le \frac{1}{2} \ln^2 d + c_2,
$$
where $c_1$ and $c_2$ are absolute constants.
Thus, the asymptotic behaviour of $p_d(n)$ admits the estimate
$$
\alpha_d \le d^{\frac{1}{2} \ln d} \alpha_2.
$$
This is better than our bound $\alpha_d\le d^{\ln d}\alpha_2$, which follows from \rf[up1], by the factor $\frac{1}{2}$ at the logarithm $\ln d$.

We note that, from our proof of Theorem~\ref{thm4} for the uniform bound, we can extract a similar estimate,
namely from \rf[eqn:newbd] we could derive 
$$
\ln \alpha_d   \le  \ln \alpha_2 + \sum_{k=3}^d \frac{\ln (k+1)}{k-1} ,
$$
which is only slightly larger than \rf[eqn:1], and will give the same asymptotic bound 
with factor $\frac{1}{2}$. 

We did not go for such improvement for two reasons: firstly, we wanted to have a uniform 
upper bound to be as simple as possible, and secondly, as we mentioned at the end of Introduction, we believe that the actual value 
of $\alpha_d$ grows much slower than $d^{c\ln d}$ (it may well be bounded by an absolute constant).
\end{remark}


\section{Universal discretization of the $L_2$ norm of trigonometric polynomials with harmonics in the lower sets}


In this section we apply our results on the cardinality of lower sets to the problem of universal sampling discretization. Currently, sampling discretization of integral norms is an actively developing area. For the background and motivation, we refer the reader to the following two recent survey papers on the topic \cite{d} and \cite{KKLT}. 

Consider the space $L_2(K)$ endowed with some probability measure on the compact
set $K \subset \R^d$.
Given an $n$-dimensional subspace $U_n \subset L_2(K) \cap C(K)$, the {\it Marcinkievicz discretization problem}
consists of finding a discrete point set $\xi = \{\xi_i\}_{i=1}^m$ such that,
for any $f \in U_n$, we have
$$
    c_1 \|f\|_2^2 \le \frac{1}{m} \sum_{i=1}^m |f(\xi_i)|^2 \le c_2 \|f\|_2^2\,,
$$
with some constants $c_1,c_2 > 0$ independent of $f$.

Now, suppose that we are given a collection of $n$-dimensional subspaces
$\U_A := \{U_n^\alpha\}_{\alpha \in A}$.
We call a point set $\xi = \{\xi_i\}_{i=1}^m$ a {\it universal discretization set} if, for any $\alpha \in A$
and for any $f \in U_n^\alpha$, we have
$$
    c_1 \|f\|_2^2 \le \frac{1}{m} \sum_{i=1}^m |f(\xi_i)|^2 \le c_2 \|f\|_2^2\,,
$$
with some constants $c_1,c_2 > 0$ independent of $f$ and $\alpha$

The following theorem provides an upper estimate for the number $m$ of points
$\{\xi_i\}_{i=1}^m$ such that, under certain
assumptions regarding the subspaces $U_n^\alpha$,
a universal discretization set for the collection $\U_A = \{U_n^\alpha\}_{\alpha \in A}$ exists. This result readily follows from minor modifications of the arguments preceding Theorem 7.4 in~\cite{d}.

\begin{theorem} \lb{thm5}
Suppose that each $U_n^\alpha \in \U_A$ is spanned by an orthonormal system $\{u_i^\alpha\}$ that satisfies
the so-called {\it Condition~E: there exists a constant $B$ such that }
$$
      \sum_{i=1}^n |u_i^\alpha|^2 \le B n.
$$
Then there exists a universal discretization set $\xi = \{\xi_i\}_{i=1}^m$ such that, for any $\alpha \in A$
and for any $f \in U_n^\alpha$, we have
$$
    c_1 \|f\|_2^2 \le \frac{1}{m} \sum_{i=1}^m |f(\xi_i)|^2 \le c_2 \|f\|_2^2\,,
$$
where 
$$
      m \le c_3 n \ln (n |A| )\,,
$$
and $c_1,c_2,c_3$ depend only on $B$.
\end{theorem}


 Our goal is to find the cardinality $m$ of a universal discretization set $\xi$ for
subspaces spanned by orthogonal systems whose indexes are from lower
sets in $\Z^d $ of size $n$.
By Theorem \ref{thm5}, we can estimate this $m$ using the results about the
cardinality of the set of all lower sets of size $n$ that we just obtained.

\begin{theorem} \lb{thm6}
For any $n$ and $d$ and any collection of subspaces $\U_A$, which are spanned by orthogonal systems whose indexes are from lower sets in $\Z^d $ of size $n$ 
and satisfy Condition~E, there is a universal discretization set  $\xi = (\xi_i)_{i=1}^m \subset K$  with the cardinality $m$ of $\xi$ satisfying the following upper bounds
\be \lb{1}
\begin{array}{l@{}l}
    a) \quad p_d(n) \le d^{n-1} &\ra m \le c n^2 \ln d,  \\
    b) \quad  p_d(n) \le \gamma_d^{n^{1-1/d}} &\ra m \le c n^{2-1/d} d^{\ln d}.
\end{array}
\ee
\end{theorem}

\proof
By Theorem \ref{thm5},
$$
       m < c_3 n \ln (n p_d(n)) = c_3 n \ln n + c_3 n \ln p_d(n)\,,
$$
and by Theorems \ref{thm3} and \ref{thm4},
$$
       a) \; \ln p_d(n) \le n \ln d, \qquad b) \; \ln p_d(n) \le \alpha_2 n^{1-1/d} d^{\ln d},
$$
so those term dominate over $\ln n$.
\qed

\medskip

Let us discuss a specific setting. Let $\{u_k(x)\}_{k=1}^\infty$ be an orthonormal system on $[0,1]$ endowed with a probability measure $\mu$. Consider $K =[0,1]^d$ endowed with the probability measure $\mu^d =\mu\times\cdots\times\mu$ and the system $\{u_\bk(\bx)\}_{\bk\in\N^d}$, where  $\bk =(k_1,\dots,k_d)$, $\bx =(x_1,\dots,x_d)$, $u_\bk(\bx):= u_{k_1}(x_1)\cdots u_{k_d}(x_d)$. Suppose that $\alpha$ stands for a lower set of the cardinality $n$. Also, assume that for each $\alpha$ Condition~E is satsified:
$$
\sum_{\bk\in\alpha}|u_\bk(\bx)|^2 \le Bn.
$$
Then Theorem \ref{thm6} can be applied to this collection.
The upper bound for $m$ provided by Theorem \ref{thm6} might not be optimal. Therefore, it makes sense to
compare our bounds \rf[1] with the bound that could be obtained by other means. All positive lower sets of size $n$
are contained in the hyperbolic cross
$$
      H_d^n := \{\k \in \N^d: \prod_{i=1}^d k_i \le n\}\,,
$$
hence, we can get a universal discretization $\xi$ of all lower sets 
just by choosing a set of points
$\xi = (\xi_i)_{i=1}^m$ that discretizes {\it all} polynomials with harmonics in the hyperbolic cross $H_d^n$.
This can be achieved with $m < c\, |H_d^n|$.
For the cardinality of $H_d^n$, we have a uniform upper bound \cite{m} $|H_d^n| <  n(1+\ln n)^{d-1}$, hence
\be \lb{H}
      m <  c\, |H_d^n| < c n(1+\ln n)^{d-1}\,.
\ee
It is clear that our estimates \rf[1] would be better than that only for large $d$ and relatively small $n$.

Let us find an approximate range of $n$ and $d$ when the estimate (a) in \rf[1] is superior to \rf[H].
For simplicity, we change the exponent on the right-hand side of \rf[H] from $d-1$ to $d$,
and ignore constants. Thus we need
\be \lb{ln}
      n^2 \ln d < n (1+\ln n)^d \lr
     \frac{(1+\ln n)^d}{n} > \ln d\,.
\ee
With $d$ fixed, consider the function $f(x) := \frac{(1+\ln x)^d}{x}$.
On $[1,\infty)$, it has one local maximum at $x = e^{d-1}$, and $f(x) \to 0$  as $x \to \infty$, and
$f(2) = \frac{(1+\ln 2)^d}{2} > \ln d$. Thus, if \rf[ln] is true for $n = n_0$, then it is true for all $2 \le n \le n_0$.

Taking $n_0 = d^d$, we check with \rf[ln]
$$
     \frac{(1 + d \ln d)^d}{d^d} > \ln d\,,
$$
which is certainly true. We may also take a bit larger $n_0 = (d\ln d)^d$.

Thus, for a universal discretization of the lower sets, we have the estimates
$$
     m < \left\{ \begin{array}{ll}
     c n^2 \ln d, &  n < d^d, \\
     c n (1+\ln n)^{d-1},  & n > d^d.
     \end{array} \right.
$$

As to the estimate (b) in \rf[1], we see that it is better than (a) when
$$
    n^2 \ln d > n^{2-1/d} d^{\ln d} \lr n > d^{d\ln d}
$$
but then, for such $n$, it would not beat \rf[H].

\begin{remark}
Theorem \ref{thm6} is based on Theorem \ref{thm5}, which in turn is a corollary of 
deep results on random matrices (see the corresponding discussion in \cite{VT159}). As a result, Theorem \ref{thm5} can be applied to any collection of subspaces of $L_2(\Omega,\mu)$ spanned by orthogonal systems satisfying Condition~E. Clearly, a uniformly bounded orthonormal system satisfies 
Condition~E.  In the very recent paper \cite{DT} the following more specific collections of subspaces of $L_2(K,\mu)$ were considered. We formulate that setting in a way convenient for our discussion. Let $\Phi_N=\{\varphi_i\}_{i=1}^N$ be a uniformly bounded orthonormal system. As a collection for the universal discretization consider
$$
\Sigma_n(\Phi_N):= \left\{f\,:\, f = \sum_{i\in G} c_i\varphi_i,\quad\text{with any}\, G\subset [1,N],\quad |G|=n\right\}.
$$
It was proved in \cite{DT} that under the above assumptions one can obtain the universal discretization for the collection $\Sigma_n(\Phi_N)$ with the following bound on the number of points
$$
m \le Cn(\log N)^2(\log(2n))^2.
$$
\end{remark}

{\bf Acknowledgement.} The authors thank the referees for their valuable comments.


 \Addresses

\end{document}